# Infinite-dimensional Log-Determinant divergences between positive definite Hilbert-Schmidt operators

Hà Quang Minh

*Istituto Italiano di Tecnologia, Via Morego 30, Genova 16163, ITALY*

**Abstract**

The current work generalizes the author's previous work on the infinite-dimensional Alpha Log-Determinant (Log-Det) divergences and Alpha-Beta Log-Det divergences, defined on the set of positive definite unitized trace class operators on a Hilbert space, to the entire Hilbert manifold of positive definite unitized Hilbert-Schmidt operators. This generalization is carried out via the introduction of the extended Hilbert-Carleman determinant for unitized Hilbert-Schmidt operators, in addition to the previously introduced extended Fredholm determinant for unitized trace class operators. The resulting parametrized family of Alpha-Beta Log-Det divergences is general and contains many divergences between positive definite unitized Hilbert-Schmidt operators as special cases, including the infinite-dimensional affine-invariant Riemannian distance and the infinite-dimensional generalization of the symmetric Stein divergence.

*Email address:* `minh.haquang@iit.it` (Hà Quang Minh)

# 1. Introduction

The current work is a continuation and generalization of the author's previous work [1], [2], which generalizes the finite-dimensional Log-Determinant divergences to the infinite-dimensional setting. We recall that for the convex cone $\text{Sym}^{++}(n)$ of symmetric, positive definite (SPD) matrices of size $n \times n$, $n \in \mathbb{N}$, the Alpha-Beta Log-Determinant (Log-Det) divergence between $A, B \in \text{Sym}^{++}(n)$ is a parametrized family of divergences defined by (see [3])

$$D^{(\alpha,\beta)}(A,B) = \frac{1}{\alpha\beta} \log \det \left[ \frac{\alpha(AB^{-1})^\beta + \beta(AB^{-1})^{-\alpha}}{\alpha + \beta} \right], \alpha > 0, \beta > 0, \quad (1)$$

along with the limiting cases $(\alpha > 0, \beta = 0), (\alpha = 0, \beta > 0)$, and $(\alpha = 0, \beta = 0)$. This family contains many distance-like functions on $\text{Sym}^{++}(n)$, including

1. The affine-invariant Riemannian distance $d_{\text{aiE}}$ [4], corresponding to

    $$D^{(0,0)}(A,B) = \frac{1}{2}d_{\text{aiE}}^2(A,B) = \frac{1}{2}||\log(B^{-1/2}AB^{-1/2})||_F^2, \quad (2)$$

    where $\log(A)$ denotes the principal logarithm of the matrix $A$ and $||\;||_F$ denotes the Frobenius norm. This is the geodesic distance associated with the so-called affine-invariant Riemannian metric [5, 6, 4, 7, 8].

2. The Alpha Log-Det divergences [9], corresponding to $D^{(\alpha,1-\alpha)}(A,B)$, with

    $$D^{(\alpha,1-\alpha)}(A,B) = \frac{1}{\alpha(1-\alpha)} \log \left[ \frac{\det[\alpha A + (1-\alpha)B]}{\det(A)^\alpha \det(B)^{1-\alpha}} \right], 0 < \alpha < 1, \quad (3)$$

    $$D^{(1,0)}(A,B) = \text{tr}(A^{-1}B - I) - \log\det(A^{-1}B), \quad (4)$$

    $$D^{(0,1)}(A,B) = \text{tr}(B^{-1}A - I) - \log\det(B^{-1}A). \quad (5)$$

    The case $\alpha = 1/2$ gives the symmetric Stein divergence (also called the Jensen-Bregman LogDet divergence), whose square root is a metric on $\text{Sym}^{++}(n)$ [10], with $D^{(1/2,1/2)}(A,B) = 4d_{\text{stein}}^2(A,B) = 4[\log\det(\frac{A+B}{2}) - \frac{1}{2}\log\det(AB)]$.

**Previous work**. In [1], we generalized the Alpha Log-Det divergences between SPD matrices [9] to the infinite-dimensional Alpha Log-Determinant divergences between positive definite unitized trace class operators on an infinite-dimensional Hilbert space. This is done via the introduction of the extended Fredholm determinant for unitized trace class operators, along with the corresponding generalization of the log-concavity of the determinant for SPD matrices to the infinite-dimensional setting. In [2], we present a formulation for the Alpha-Beta Log-Det divergences between positive definite unitized trace class operators, generalizing the Alpha-Beta Log-Det divergences between SPD matrices as defined by Eq.(1). In both [1] and [2], for the divergences between reproducing kernel Hilbert spaces (RKHS) covariance operators, we obtain closed form formulas for the Alpha-Beta Log-Det divergences via the corresponding Gram matrices.

**Contributions of this work**. The current work is a continuation and generalization of [1] and [2]. In particular, we generalize the Alpha-Beta Log-Det divergences in [2] to the entire Hilbert manifold of positive definite unitized Hilbert-Schmidt operators on an infinite-dimensional Hilbert space. This is done by the introduction of the extended Hilbert-Carleman determinant for unitized Hilbert-Schmidt operators, in addition to the extended Fredholm determinant for unitized trace class operators employed in [1] and [2]. As in the finite-dimensional setting [3] and in [1], [2], the resulting family of divergences is general and admits as special cases many metrics and distance-like functions between positive definite unitized Hilbert-Schmidt operators, including the infinite-dimensional affine-invariant Riemannian distance in [11].

**Comparison with the formulations in [1] and [2]**. While the mathematical formulation presented in the current work, for Hilbert-Schmidt operators, is more general than the formulations in [1] and [2], which are for trace class operators, it should *not* be considered as a substitute for them. Many results in [1] and [2], especially those involving covariance operators, require explicitly the trace class assumption.

## 2. Positive definite unitized trace class and Hilbert-Schmidt operators

Throughout the paper, we assume that $\mathcal{H}$ is a real separable Hilbert space, with $\dim(\mathcal{H}) = \infty$, unless explicitly stated otherwise. Let $\mathcal{L}(\mathcal{H})$ be the Banach space of bounded linear operators on $\mathcal{H}$, with operator norm $||\ ||$. Let $\mathrm{Sym}(\mathcal{H}) \subset \mathcal{L}(\mathcal{H})$ denote the subspace of bounded, self-adjoint operators on $\mathcal{H}$. Let $\mathrm{Sym}^+(\mathcal{H}) \subset \mathrm{Sym}(\mathcal{H})$ denote the set of self-adjoint, *positive* operators on $\mathcal{H}$, that is $A \in \mathrm{Sym}^+(\mathcal{H}) \iff$ $\langle x, Ax \rangle \geq 0 \ \forall x \in \mathcal{H}$. Let $\text{Sym}^{++}(\mathcal{H}) \subset \text{Sym}^+(\mathcal{H})$ denote the set of self-adjoint, *strictly positive* operators on $\mathcal{H}$, that is $A \in \text{Sym}^{++}(\mathcal{H}) \iff \langle x, Ax \rangle > 0 \ \forall x \in \mathcal{H}, x \neq 0$, or equivalently, $\ker(A) = \{0\}$.

Most importantly, we consider the set $\mathbb{P}(\mathcal{H}) \subset \text{Sym}^{++}(\mathcal{H})$ of self-adjoint, bounded, *positive definite* operators on $\mathcal{H}$, which is defined by

$$A \in \mathbb{P}(\mathcal{H}) \iff A = A^*, \exists M_A > 0 \text{ such that } \langle x, Ax \rangle \geq M_A ||x||^2 \quad \forall x \in \mathcal{H}.$$

We use the notation $A > 0 \iff A \in \mathbb{P}(\mathcal{H})$.

In the following, let $\mathscr{C}_p(\mathcal{H})$ denote the set of $p$th Schatten class operators on $\mathcal{H}$ (see e.g. [12]), under the norm $|| \ ||_p$, $1 \leq p \leq \infty$, which is defined by

$$\mathscr{C}_p(\mathcal{H}) = \{A \in \mathcal{L}(\mathcal{H}) \ : \ ||A||_p = (\text{tr}|A|^p)^{1/p} < \infty\}, \tag{6}$$

where $|A| = (A^*A)^{1/2}$.

The cases we consider in this work are: (i) the space $\mathscr{C}_1(\mathcal{H})$ of trace class operators on $\mathcal{H}$, which we also denote by $\text{Tr}(\mathcal{H})$, and (ii) the space $\mathscr{C}_2(\mathcal{H})$ of Hilbert-Schmidt operators on $\mathcal{H}$, which we also denote by $\text{HS}(\mathcal{H})$.

**Extended (unitized) trace class operators**. In [1], we define the set of extended (or unitized) trace class operators on $\mathcal{H}$ to be

$$\text{Tr}_X(\mathcal{H}) = \{A + \gamma I \ : \ A \in \text{Tr}(\mathcal{H}), \gamma \in \mathbb{R}\}.$$

The set $\text{Tr}_X(\mathcal{H})$ becomes a Banach algebra under the *extended trace class norm*

$$||A + \gamma I||_{\text{tr}_X} = ||A||_{\text{tr}} + |\gamma| = \text{tr}|A| + |\gamma|.$$

For $(A + \gamma I) \in \text{Tr}_X(\mathcal{H})$, its *extended trace* is defined to be

$$\text{tr}_X(A + \gamma I) = \text{tr}(A) + \gamma.$$

By this definition $\text{tr}_X(I) = 1$, in contrast to standard trace definition, according to which $\text{tr}(I) = \infty$.

**Extended (unitized) Hilbert-Schmidt operators**. In [11], the author considered the following set of extended (unitized) Hilbert-Schmidt operators

$$\text{HS}_X(\mathcal{H}) = \{A + \gamma I \ : \ A \in \text{HS}(\mathcal{H}), \gamma \in \mathbb{R}\}. \tag{7}$$

The set $\text{HS}_X(\mathcal{H})$ can be equipped with the *extended Hilbert-Schmidt inner product* $\langle\ ,\ \rangle_{\text{eHS}}$, defined by

$$\langle A + \gamma I, B + \mu I\rangle_{\text{eHS}} = \langle A, B\rangle_{\text{HS}} + \gamma\mu = \text{tr}(A^*B) + \gamma\mu.$$

along with the associated *extended Hilbert-Schmidt norm*

$$||A + \gamma I||^2_{\text{eHS}} = ||A||^2_{\text{HS}} + \gamma^2 = \text{tr}(A^*A) + \gamma^2. \tag{8}$$

Under the inner product $\langle\ ,\ \rangle_{\text{eHS}}$, the Hilbert-Schmidt operators are orthogonal to the scalar operators. Under the norm $||\ ||_{\text{eHS}}$, $||I||_{\text{eHS}} = 1$, in contrast to the standard Hilbert-Schmidt norm, according to which $||I||_{\text{HS}} = \infty$.

**Positive definite unitized trace class and Hilbert-Schmidt operators**. The set of positive definite unitized trace class operators $\mathscr{PC}_1(\mathcal{H}) \subset \text{Tr}_X(\mathcal{H})$ is defined to be the intersection

$$\mathscr{PC}_1(\mathcal{H}) = \text{Tr}_X(\mathcal{H}) \cap \mathbb{P}(\mathcal{H}) = \{A + \gamma I > 0 \ :\ A^* = A,\ A \in \text{Tr}(\mathcal{H}),\ \gamma \in \mathbb{R}\}. \tag{9}$$

The set of positive definite unitized Hilbert-Schmidt operators $\mathscr{PC}_2(\mathcal{H}) \subset \text{HS}_X(\mathcal{H})$ is defined to be the intersection

$$\mathscr{PC}_2(\mathcal{H}) = \text{HS}_X(\mathcal{H}) \cap \mathbb{P}(\mathcal{H}) = \{A + \gamma I > 0 \ :\ A = A^*, A \in \text{HS}(\mathcal{H}), \gamma \in \mathbb{R}\}. \tag{10}$$

*Remark* 1. In [1] and [2], we use the notations $\text{PTr}(\mathcal{H})$ and $\Sigma(\mathcal{H})$ to denote $\mathscr{PC}_1(\mathcal{H})$ and $\mathscr{PC}_2(\mathcal{H})$, respectively. In the following, we refer to elements of $\mathscr{PC}_1(\mathcal{H})$ and $\mathscr{PC}_2(\mathcal{H})$ as *positive definite trace class operators* and *positive definite Hilbert-Schmidt operators*, respectively.

In [11], it is shown that the set $\mathscr{PC}_2(\mathcal{H})$ assumes the structure of an infinite-dimensional Hilbert manifold and can be equipped with the following Riemannian metric. For each $P \in \mathscr{PC}_2(\mathcal{H})$, on the tangent space $T_P(\mathscr{PC}_2(\mathcal{H})) \cong \mathcal{H}_\mathbb{R} = \{A + \gamma I\ :\ A = A^*, A \in \text{HS}(\mathcal{H}), \gamma \in \mathbb{R}\}$, we define the following inner product

$$\langle A + \gamma I, B + \mu I\rangle_P = \langle P^{-1/2}(A + \gamma I)P^{-1/2}, P^{-1/2}(B + \mu I)P^{-1/2}\rangle_{\text{eHS}}.$$

The Riemannian metric given by $\langle\ ,\ \rangle_P$ then makes $\mathscr{PC}_2(\mathcal{H})$ an infinite-dimensional Riemannian manifold. Under this Riemannian metric, the geodesic distance between $(A+\gamma I), (B+\mu I) \in \mathscr{PC}_2(\mathcal{H})$ is given by

$$d_{\text{aiHS}}[(A+\gamma I),(B+\mu I)] = ||\log[(B+\mu I)^{-1/2}(A+\gamma I)(B+\mu I)^{-1/2}]||_{\text{eHS}}. \tag{11}$$

**Aim of this work**. In [1], we introduce a parametrized family of divergences, called *Log-Determinant divergences*, between operators in $\mathscr{PC}_1(\mathcal{H})$. In [2], we generalize these to the *Alpha-Beta Log-Determinant divergences* on $\mathscr{PC}_1(\mathcal{H})$, which include the distance $d_{\text{aiHS}}$ as a special case. However, these divergences are defined specifically on $\mathscr{PC}_1(\mathcal{H})$. In the case $\dim(\mathcal{H}) = \infty$, the set $\mathscr{PC}_1(\mathcal{H})$ of positive definite trace class operators on $\mathcal{H}$ is a *strict subset* of the set of positive definite Hilbert-Schmidt operators $\mathscr{PC}_2(\mathcal{H})$. In this work, we generalize the divergences in [1] and [2] to all of $\mathscr{PC}_2(\mathcal{H})$.

### 3. Functions of positive definite unitized Hilbert-Schmidt operators

We first discuss several important functions on $\mathscr{PC}_2(\mathcal{H})$, namely the exponential, logarithm, and power functions.

**Exponential and logarithm functions**. Consider the exponential function $\exp : \mathcal{L}(\mathcal{H}) \to \mathcal{L}(\mathcal{H})$ defined by

$$\exp(A) = \sum_{k=0}^{\infty} \frac{A^k}{k!}. \tag{12}$$

In [11], it is shown that the map $\exp : \text{Sym}(\mathcal{H}) \cap \text{HS}_X(\mathcal{H}) \to \mathscr{PC}_2(\mathcal{H})$ and its inverse function $\log = \exp^{-1} : \mathscr{PC}_2(\mathcal{H}) \to \text{Sym}(\mathcal{H}) \cap \text{HS}_X(\mathcal{H})$ are diffeomorphisms. Here, for any $(A+\gamma I) \in \mathscr{PC}_2(\mathcal{H})$, $\log(A+\gamma I)$ is defined via the spectral decomposition of $A$ as follows. Let $\{\lambda_k\}_{k=1}^{\infty}$ be the eigenvalues of $A$ with corresponding orthonormal eigenvectors $\{\phi_k\}_{k=1}^{\infty}$. Then

$$A = \sum_{k=1}^{\infty} \lambda_k \phi_k \otimes \phi_k, \quad \log(A+\gamma I) = \sum_{k=1}^{\infty} \log(\lambda_k + \gamma)\phi_k \otimes \phi_k, \tag{13}$$

where $\phi_k \otimes \phi_k : \mathcal{H} \rightarrow \mathcal{H}$ is a rank-one operator defined by $(\phi_k \otimes \phi_k)w = \langle \phi_k, w\rangle \phi_k$ $\forall w \in \mathcal{H}$. Since $\log(A + \gamma I) \in \text{Sym}(\mathcal{H}) \cap \text{HS}_X(\mathcal{H})$, it has the form

$$\log(A + \gamma I) = A_1 + \gamma_1 I, \quad A_1 \in \text{Sym}(\mathcal{H}) \cap \text{HS}(\mathcal{H}), \gamma_1 \in \mathbb{R}.$$

**Power functions**. Given the exponential and logarithm functions, for any $\alpha \in \mathbb{R}$, the power function $(A + \gamma I)^\alpha$, for $(A + \gamma I) \in \mathscr{PC}_2(\mathcal{H})$, is then well-defined via the following expression

$$(A + \gamma I)^\alpha = \exp[\alpha \log(A + \gamma I)] \in \mathscr{PC}_2(\mathcal{H}).$$

Furthermore, for any two operators $(A + \gamma I), (B + \mu I) \in \mathscr{PC}_2(\mathcal{H})$, we show that

$$\log[(A + \gamma I)(B + \mu I)^{-1}], \quad [(A + \gamma I)(B + \mu I)^{-1}]^\alpha, \alpha \in \mathbb{R} \quad (14)$$

are all well-defined and are elements of $\text{HS}_X(\mathcal{H})$ (though not necessarily of $\text{Sym}(\mathcal{H})$). To this end, let $B \in \mathcal{L}(\mathcal{H})$ be any invertible operator, then for any $A \in \mathcal{L}(\mathcal{H})$, we have

$$\exp(BAB^{-1}) = \sum_{j=0}^\infty \frac{(BAB^{-1})^j}{j!} = B\left(\sum_{j=0}^\infty \frac{A^j}{j!}\right)B^{-1} = B\exp(A)B^{-1}.$$

Thus for $(A + \gamma I) \in \mathscr{PC}_2(\mathcal{H})$, the logarithm of $B(A + \gamma I)B^{-1} = BAB^{-1} + \gamma I \in \text{HS}_X(\mathcal{H})$ is also well-defined and is given by

$$\log[B(A + \gamma I)B^{-1}] = B\log(A + \gamma I)B^{-1}$$
$$= B(A_1 + \gamma_1 I)B^{-1} = BA_1 B^{-1} + \gamma_1 I \in \text{HS}_X(\mathcal{H}). \quad (15)$$

Using Eq. (15), we obtain the following results.

**Proposition 1.** *Let* $(A+\gamma I), (B+\mu I) \in \mathscr{PC}_2(\mathcal{H})$. *Let* $\Lambda + \frac{\gamma}{\mu}I = (B+\mu I)^{-1/2}(A+\gamma I)(B+\mu I)^{-1/2}$. *Then*

1. *The logarithm function* $\log[(A + \gamma I)(B + \mu I)^{-1}] \in \text{HS}_X(\mathcal{H})$ *is well-defined and is given by*

$$\log[(A + \gamma I)(B + \mu I)^{-1}] = (B + \mu I)^{1/2} \log\left(\Lambda + \frac{\gamma}{\mu}I\right)(B + \mu I)^{-1/2}. \quad (16)$$

2. *For any $\alpha \in \mathbb{R}$, the power function $[(A + \gamma I)(B + \mu I)^{-1}]^{\alpha} \in \mathrm{HS}_X(\mathcal{H})$ is well-defined and is given by*

$$[(A + \gamma I)(B + \mu I)^{-1}]^{\alpha} = (B + \mu I)^{1/2} \left(\Lambda + \frac{\gamma}{\mu}I\right)^{\alpha} (B + \mu I)^{-1/2}. \quad (17)$$

**4. The extended Hilbert-Carleman determinant**

The key concept for defining Log-Determinant divergences between operators is determinant. We recall that for $A \in \mathrm{Tr}(\mathcal{H})$, the Fredholm determinant $\det(I + A)$ is (see e.g. [13])

$$\det(I + A) = \prod_{k=1}^{\infty}(1 + \lambda_k), \quad (18)$$

where $\{\lambda_k\}_{k=1}^{\infty}$ are the eigenvalues of $A$. To define Log-Determinant divergences between positive definite trace class operators in $\mathscr{PC}_1(\mathcal{H})$, in [1], we generalize the Fredholm determinant to the *extended Fredholm determinant* of extended trace class operators. For $(A + \gamma I) \in \mathrm{Tr}_X(\mathcal{H})$, $\gamma \neq 0$, its extended Fredholm determinant is defined to be, assuming that $\dim(\mathcal{H}) = \infty$,

$$\det{}_X(A + \gamma I) = \frac{1}{\gamma}\det\left(\frac{A}{\gamma} + I\right),$$

where the determinant on the right hand side is the Fredholm determinant (we refer to [1] for the derivation leading to this definition). For $\gamma = 1$, we recover the Fredholm determinant. In the case $\dim(\mathcal{H}) < \infty$, we define $\det{}_X(A + \gamma I) = \det(A + \gamma I)$, the standard matrix determinant.

The extended Fredholm determinant continues to play a key role in the current work, but it is not sufficient for dealing with positive definite Hilbert-Schmidt operators in $\mathscr{PC}_2(\mathcal{H})$. In order to do so, we introduce the concept of *extended Hilbert-Carleman determinant*.

We first recall the concept of the Hilbert-Carleman determinant for operators of the form $I + A$, where $A$ is a Hilbert-Schmidt operator (see e.g. [13] for a comprehensive treatment). Following [13], for any bounded operator $A \in \mathcal{L}(\mathcal{H})$, consider the operator

$$R_n(A) = \left[(I + A)\exp\left(\sum_{k=1}^{n-1}\frac{(-A)^k}{k}\right)\right] - I. \quad (19)$$

If $A \in \mathscr{C}_n(\mathcal{H})$, then $R_n(A) \in \mathscr{C}_1(\mathcal{H})$. Thus the following quantity is well-defined

$$\det{}_n(I+A) = \det(I + R_n(A)). \tag{20}$$

In particular, for $n = 1$, we obtain $R_1(A) = A$ and thus

$$\det{}_1(I+A) = \det(I+A). \tag{21}$$

For $n = 2$, we have $R_2(A) = (I+A)\exp(-A) - I$ and thus

$$\det{}_2(I+A) = \det[(I+A)\exp(-A)]. \tag{22}$$

This is called the *Hilbert-Carleman determinant* of $I + A$. In particular, for $A \in \mathrm{Tr}(\mathcal{H}) = \mathscr{C}_1(\mathcal{H})$, we have

$$\det{}_2(I+A) = \det(I+A)\exp(-\mathrm{tr}(A)), \tag{23}$$

$$\log\det{}_2(I+A) = \log\det(I+A) - \mathrm{tr}(A). \tag{24}$$

The function $\det_2(I+A)$ is continuous in the Hilbert-Schmidt norm, so that

$$\lim_{k\to\infty} ||A_k - A||_{\mathrm{HS}} = 0 \Rightarrow \lim_{k\to\infty} \det{}_2(I+A_k) = \det{}_2(I+A). \tag{25}$$

We first have the following result.

**Lemma 1.** *Assume that $A \in \mathrm{Sym}(\mathcal{H}) \cap \mathrm{HS}(\mathcal{H})$ such that $I + A > 0$. Let $\{\lambda_k\}_{k=1}^{\infty}$ be the eigenvalues of $A$. Then*

$$\log\det{}_2(I+A) = \sum_{k=1}^{\infty}[\log(1+\lambda_k) - \lambda_k] \tag{26}$$

*is well-defined and finite. Furthermore,*

$$\log\det{}_2(I+A) \leq 0, \tag{27}$$

*with equality if and only if $A = 0$.*

The Hilbert-Carleman determinant $\det_2$ is defined for operators of the form $A + I$, $A \in \mathrm{HS}(\mathcal{H})$, but not for operators of the form $A + \gamma I$, $\gamma > 0, \gamma \neq 1$. In the following, we generalize $\det_2$ to handle these operators. We first have the following generalization of the function $R_2(A) = (I+A)\exp(-A) - I$ above.

**Lemma 2.** *Assume that $(A + \gamma I) \in \mathrm{HS}_X(\mathcal{H})$, $\gamma \neq 0$. Define*

$$R_{2,\gamma}(A) = (A + \gamma I) \exp(-A/\gamma) - \gamma I. \tag{28}$$

*Then $R_{2,\gamma}(A) \in \mathrm{Tr}(\mathcal{H})$ and hence $R_{2,\gamma}(A) + \gamma I = (A+\gamma I)\exp(-A/\gamma) \in \mathrm{Tr}_X(\mathcal{H})$. This also implies that the infinite product*

$$\prod_{k=1}^{\infty} [(\lambda_k + \gamma)\exp(-\lambda_k/\gamma) - \gamma + 1] \tag{29}$$

*converges to a finite value, where $\{\lambda_k\}_{k=1}^{\infty}$ are the eigenvalues of $A$.*

In particular, for $\gamma = 1$, we have $R_{2,1}(A) = R_2(A)$. Motivated by Lemma 2 and the definition of $\det_2$, we arrive at the following generalization of $\det_2$.

**Definition 1** (**Extended Hilbert-Carleman determinant**)**.** *For $(A+\gamma I) \in \mathrm{HS}_X(\mathcal{H})$, $\gamma \neq 0$, its extended Hilbert-Carleman determinant is defined to be*

$$\det{}_{2\mathrm{X}}(A + \gamma I) = \det{}_{\mathrm{X}}[R_{2,\gamma}(A) + \gamma I] = \det{}_{\mathrm{X}}[(A+\gamma I)\exp(-A/\gamma)]. \tag{30}$$

If $\gamma = 1$, then we recover the Hilbert-Carleman determinant

$$\det{}_{2\mathrm{X}}(A + I) = \det[(A + I)\exp(-A)] = \det{}_2(A + I). \tag{31}$$

If $(A + \gamma I) \in \mathrm{Tr}_X(\mathcal{H})$, then

$$\det{}_{2\mathrm{X}}(A + \gamma I) = \det{}_{\mathrm{X}}(A + \gamma I)\exp(-\mathrm{tr}(A)/\gamma). \tag{32}$$

The following are the some of the properties of $\det_{2\mathrm{X}}$ which we employ later on.

**Lemma 3** (**Factorization Rule**)**.**

$$\det{}_{2\mathrm{X}}(A + \gamma I) = \gamma \det{}_2\left(\frac{A}{\gamma} + I\right). \tag{33}$$

If $(A + \gamma I) \in \mathrm{Tr}_X(\mathcal{H})$, $\gamma \neq 0$, then Lemma 5 in [1] states that for any invertible operator $C \in \mathcal{L}(\mathcal{H})$, we have

$$\det{}_{\mathrm{X}}[C(A + \gamma I)C^{-1}] = \det{}_{\mathrm{X}}(A + \gamma I). \tag{34}$$

This property generalizes for $\det_{2\mathrm{X}}$, with $(A + \gamma I) \in \mathrm{HS}_X(\mathcal{H})$, as follows.

**Lemma 4 (Similarity Invariant).** *Let $(A + \gamma I) \in \mathrm{HS}_X(\mathcal{H})$, $\gamma \neq 0$. Let $C \in \mathcal{L}(\mathcal{H})$ be invertible. Then*

$$\mathrm{det}_{2\mathrm{X}}[C(A + \gamma I)C^{-1}] = \mathrm{det}_{2\mathrm{X}}(A + \gamma I). \tag{35}$$

For $(A + \gamma I), (B + \mu I) \in \mathrm{Tr}_X(\mathcal{H})$, we show in Proposition 4 in [1] that the product rule for determinants holds, that is $\mathrm{det}_\mathrm{X}[(A + \gamma I)(B + \mu I)] = \mathrm{det}_\mathrm{X}(A + \gamma I)\mathrm{det}_\mathrm{X}(B + \mu I)$. For $\mathrm{det}_2$ and $\mathrm{det}_{2\mathrm{X}}$ and $(A + \gamma I), (B + \mu I) \in \mathrm{HS}_X(\mathcal{H})$, this is no longer true in general. However, if $(A + \gamma I), (B + \mu I) \in \mathscr{PC}_2(\mathcal{H})$, or if $(A + \gamma I), (B + \mu I) \in \mathrm{Tr}_X(\mathcal{H})$, then we still have commutativity, that is $\mathrm{det}_{2\mathrm{X}}[(A + \gamma I)(B + \mu I)] = \mathrm{det}_{2\mathrm{X}}[(B + \mu I)(A + \gamma I)]$, as follows.

**Lemma 5 (Commutativity).** *Assume that $(A + \gamma I), (B + \mu I) \in \mathscr{PC}_2(\mathcal{H})$. Then*

$$\mathrm{det}_{2\mathrm{X}}[(A + \gamma I)(B + \mu I)] = \mathrm{det}_{2\mathrm{X}}[(A + \gamma I)^{1/2}(B + \mu I)(A + \gamma I)^{1/2}] \tag{36}$$

$$= \mathrm{det}_{2\mathrm{X}}[(B + \mu I)(A + \gamma I)] \tag{37}$$

$$= \mathrm{det}_{2\mathrm{X}}[(B + \mu I)^{1/2}(A + \gamma I)(B + \mu I)^{1/2}]. \tag{38}$$

*If $(A + \gamma I), (B + \mu I) \in \mathrm{Tr}_X(\mathcal{H})$, $\gamma \neq 0, \mu \neq 0$, then*

$$\mathrm{det}_{2\mathrm{X}}[(A + \gamma I)(B + \mu I)] = \mathrm{det}_{2\mathrm{X}}[(B + \mu I)(A + \gamma I)]. \tag{39}$$

An immediate consequence of Lemma 5 is the following.

**Corollary 1 (Cyclic Property).** *Assume that $(A + \gamma I), (B + \mu I), (C + \nu I) \in \mathscr{PC}_2(\mathcal{H})$, or $(A + \gamma I), (B + \mu I), (C + \nu I) \in \mathrm{Tr}_X(\mathcal{H})$. Then*

$$\mathrm{det}_{2\mathrm{X}}[(A + \gamma I)(B + \mu I)(C + \nu I)] = \mathrm{det}_{2\mathrm{X}}[(C + \nu I)(A + \gamma I)(B + \mu I)] \tag{40}$$

$$= \mathrm{det}_{2\mathrm{X}}[(B + \mu I)(C + \nu I)(A + \gamma I)]. \tag{41}$$

For the following properties, we assume explicitly that $(A + \gamma I), (B + \mu I) \in \mathscr{PC}_2(\mathcal{H})$, that is $(A + \gamma I) > 0, (B + \mu I) > 0$ and $A, B \in \mathrm{HS}(\mathcal{H})$. These properties are utilized in the formulation of the Log-Determinant divergences in Section 5.

**Lemma 6.** *Let* $(A + \gamma I), (B + \mu I) \in \mathscr{PC}_2(\mathcal{H})$. *Let* $\Lambda + \nu I = (B + \mu I)^{-1/2}(A + \gamma I)(B + \mu I)^{-1/2}$, *where* $\Lambda \in \mathrm{HS}(\mathcal{H})$ *and* $\nu = \frac{\gamma}{\mu}$. *Then for any* $\alpha \in \mathbb{R}$,

$$\det_{2\mathrm{X}}([(A + \gamma I)(B + \mu I)^{-1}]^\alpha) = \det_{2\mathrm{X}}[(\Lambda + \nu I)^\alpha]$$
$$= \det_{2\mathrm{X}}([(B + \mu I)^{-1}(A + \gamma I)]^\alpha). \tag{42}$$

**Lemma 7.** *Let* $(I + A) \in \mathscr{PC}_2(\mathcal{H})$. *Let* $\alpha \in \mathbb{R}$ *be arbitrary. Then* $(I + A)^\alpha - I \in \mathrm{Sym}(\mathcal{H}) \cap \mathrm{HS}(\mathcal{H})$. *Let* $\{\lambda_k\}_{k=1}^\infty$ *be the eigenvalues of A. Then the following quantities converge to finite values:*

$$\det_2[(I + A)^\alpha] = \prod_{k=1}^\infty (1 + \lambda_k)^\alpha \exp[1 - (1 + \lambda_k)^\alpha], \tag{43}$$

$$\log \det_2[(I + A)^\alpha] = \sum_{k=1}^\infty [\alpha \log(1 + \lambda_k) + 1 - (1 + \lambda_k)^\alpha]. \tag{44}$$

*Furthermore,*

$$\log \det_2[(I + A)^\alpha] \leq 0, \tag{45}$$

*with equality if and only if* $A = 0$. *For* $(I + A) \in \mathscr{PC}_1(\mathcal{H})$,

$$\det_2[(I + A)^\alpha] = \det[(I + A)^\alpha] \exp(-\mathrm{tr}[(I + A)^\alpha - I]), \tag{46}$$

$$\log \det_2[(I + A)^\alpha] = \alpha \log \det(I + A) - \mathrm{tr}[(I + A)^\alpha - I]. \tag{47}$$

**Lemma 8.** *Let* $(A + \gamma I) \in \mathscr{PC}_2(\mathcal{H})$. *Let* $\alpha \in \mathbb{R}$ *be arbitrary. Then* $(A + \gamma I)^\alpha - \gamma^\alpha I \in \mathrm{Sym}(\mathcal{H}) \cap \mathrm{HS}(\mathcal{H})$. *Let* $\{\lambda_k\}_{k=1}^\infty$ *be the eigenvalues of A. Then the following quantities converge to finite values:*

$$\det_{2\mathrm{X}}[(A + \gamma I)^\alpha] = \gamma^\alpha \prod_{k=1}^\infty \left(1 + \frac{\lambda_k}{\gamma}\right)^\alpha \exp\left[1 - \left(1 + \frac{\lambda_k}{\gamma}\right)^\alpha\right], \tag{48}$$

$$\log \det_{2\mathrm{X}}[(A + \gamma I)^\alpha] = \alpha \log \gamma + \sum_{k=1}^\infty \left[\alpha \log\left(1 + \frac{\lambda_k}{\gamma}\right) + 1 - \left(1 + \frac{\lambda_k}{\gamma}\right)^\alpha\right]. \tag{49}$$

*For* $(A + \gamma I) \in \mathscr{PC}_1(\mathcal{H})$,

$$\text{det}_{2X}[(A + \gamma I)^\alpha] = \gamma^\alpha \det\left[\left(\frac{A}{\gamma} + I\right)^\alpha\right] \exp\left(-\text{tr}\left[\left(\frac{A}{\gamma} + I\right)^\alpha - I\right]\right), \tag{50}$$

$$\log \text{det}_{2X}[(A + \gamma I)^\alpha] = \alpha \log \gamma + \alpha \log \det\left(\frac{A}{\gamma} + I\right) - \text{tr}\left[\left(\frac{A}{\gamma} + I\right)^\alpha - I\right]$$

$$= \alpha \log \text{det}_X(A + \gamma I) - \text{tr}\left[\left(\frac{A}{\gamma} + I\right)^\alpha - I\right]. \tag{51}$$

## 5. Infinite-dimensional Log-Determinant divergences between positive definite Hilbert-Schmidt operators

In [2], we define the Log-Determinant divergences between two positive definite trace class operators $(A + \gamma I), (B + \mu I) \in \mathscr{PC}_1(\mathcal{H})$ as follows

$$D_r^{(\alpha, \beta)}[(A + \gamma I), (B + \mu I)]$$
$$= \frac{1}{\alpha\beta} \log\left[\left(\frac{\gamma}{\mu}\right)^{r(\delta - \frac{\alpha}{\alpha+\beta})} \text{det}_X\left(\frac{\alpha(\Lambda + \frac{\gamma}{\mu}I)^{r(1-\delta)} + \beta(\Lambda + \frac{\gamma}{\mu}I)^{-r\delta}}{\alpha + \beta}\right)\right]. \tag{52}$$

Here $\alpha > 0, \beta > 0, r \neq 0$ are fixed, $\Lambda + \frac{\gamma}{\mu}I = (B + \mu I)^{-1/2}(A + \gamma I)(B + \mu I)^{-1/2}$, and $\delta = \frac{\alpha \gamma^r}{\alpha \gamma^r + \beta \mu^r}$. This definition is motivated by the infinite-dimensional generalizations of Ky Fan's inequality [14] on the log-concavity of the determinant of SPD matrices, as stated for $\text{det}_X$ in Theorem 1 in [1] and Theorem 5 in [2].

In the following, we show that the definition given in Eq. (52) is valid in the more general case $(A + \gamma I), (B + \mu I) \in \mathscr{PC}_2(\mathcal{H})$. We first have the following results.

**Proposition 2.** *Assume that* $(A + \gamma I), (B + \mu I) \in \mathscr{PC}_2(\mathcal{H})$. *Let* $\alpha > 0, \beta > 0$ *be fixed. Let* $p, q \in \mathbb{R}$ *be such that* $p\alpha(\gamma/\mu)^p = q\beta(\gamma/\mu)^{-q}$. *Then for* $\Lambda + \frac{\gamma}{\mu}I = (B + \mu I)^{-1/2}(A + \gamma I)(B + \mu I)^{-1/2}$,

$$\frac{\alpha(\Lambda + \frac{\gamma}{\mu}I)^p + \beta(\Lambda + \frac{\gamma}{\mu}I)^{-q}}{\alpha + \beta} \in \mathscr{PC}_1(\mathcal{H}). \tag{53}$$

**Proposition 3.** *Assume that* $(A + \gamma I), (B + \mu I) \in \mathscr{PC}_2(\mathcal{H})$. *Let* $\Lambda + \frac{\gamma}{\mu}I = (B + \mu I)^{-1/2}(A + \gamma I)(B + \mu I)^{-1/2}$. *Let* $\alpha > 0, \beta > 0$ *be fixed. Let* $p, q \in \mathbb{R}$ *be such that* $p\alpha(\gamma/\mu)^p = q\beta(\gamma/\mu)^{-q}$. *Then*

$$\frac{\alpha[(A + \gamma I)(B + \mu I)^{-1}]^p + \beta[(A + \gamma I)(B + \mu I)^{-1}]^{-q}}{\alpha + \beta} \in \text{Tr}_X(\mathcal{H}). \tag{54}$$

*Furthermore,*

$$\det_{\mathrm{X}}\left[\frac{\alpha[(A+\gamma I)(B+\mu I)^{-1}]^p + \beta[(A+\gamma I)(B+\mu I)^{-1}]^{-q}}{\alpha+\beta}\right]$$
$$=\det_{\mathrm{X}}\left[\frac{\alpha(\Lambda+\frac{\gamma}{\mu}I)^p + \beta(\Lambda+\frac{\gamma}{\mu}I)^{-q}}{\alpha+\beta}\right]. \quad (55)$$

Motivated by Eq. (52) and Propositions 2 and 3, the following is our definition of the Alpha-Beta Log-Determinant divergences on $\mathscr{PC}_2(\mathcal{H})$.

**Definition 2** (**Alpha-Beta Log-Determinant divergences between positive definite Hilbert-Schmidt operators**). *Assume that* $\dim(\mathcal{H}) = \infty$. *Let* $\alpha > 0$, $\beta > 0$ *be fixed. Let* $r \in \mathbb{R}$, $r \neq 0$ *be fixed. For* $(A+\gamma I), (B+\mu I) \in \mathscr{PC}_2(\mathcal{H})$, *the* $(\alpha, \beta)$-*Log-Det divergence* $D_r^{(\alpha,\beta)}[(A+\gamma I), (B+\mu I)]$ *is defined to be*

$$D_r^{(\alpha,\beta)}[(A+\gamma I),(B+\mu I)]$$
$$= \frac{1}{\alpha\beta}\log\left[\left(\frac{\gamma}{\mu}\right)^{r(\delta-\frac{\alpha}{\alpha+\beta})}\det_{\mathrm{X}}\left(\frac{\alpha(\Lambda+\frac{\gamma}{\mu}I)^{r(1-\delta)}+\beta(\Lambda+\frac{\gamma}{\mu}I)^{-r\delta}}{\alpha+\beta}\right)\right], \quad (56)$$

*where* $\Lambda + \frac{\gamma}{\mu}I = (B+\mu I)^{-1/2}(A+\gamma I)(B+\mu I)^{-1/2}$, $\delta = \frac{\alpha\gamma^r}{\alpha\gamma^r + \beta\mu^r}$. *Equivalently,*

$$D_r^{(\alpha,\beta)}[(A+\gamma I),(B+\mu I)]$$
$$= \frac{1}{\alpha\beta}\log\left[\left(\frac{\gamma}{\mu}\right)^{r(\delta-\frac{\alpha}{\alpha+\beta})}\det_{\mathrm{X}}\left(\frac{\alpha(Z+\frac{\gamma}{\mu}I)^{r(1-\delta)}+\beta(Z+\frac{\gamma}{\mu}I)^{-r\delta}}{\alpha+\beta}\right)\right], \quad (57)$$

*where* $Z + \frac{\gamma}{\mu}I = (A+\gamma I)(B+\mu I)^{-1}$.

While Definition 2 is stated using the extended Fredholm determinant $\det_{\mathrm{X}}$, the limiting cases $(\alpha > 0, \beta = 0)$ and $(\alpha = 0, \beta > 0)$ both require the concept of the extended Hilbert-Carleman determinant $\det_{2\mathrm{X}}$.

**Theorem 1** (**Limiting case** $\alpha > 0, \beta \to 0$). *Let* $\alpha > 0$ *be fixed. Assume that* $r = r(\beta)$ *is smooth, with* $r(0) = r(\beta = 0)$. *Then*

$$\lim_{\beta \to 0} D_r^{(\alpha,\beta)}[(A+\gamma I),(B+\mu I)] = \frac{1}{\alpha^2}\left[\left(\frac{\mu}{\gamma}\right)^{r(0)}-1\right]\left(1+r(0)\log\frac{\mu}{\gamma}\right) \quad (58)$$
$$- \frac{1}{\alpha^2}\left(\frac{\mu}{\gamma}\right)^{r(0)}\log\det_{2\mathrm{X}}([(A+\gamma I)^{-1}(B+\mu I)]^{r(0)}).$$

**Theorem 2** (**Limiting case** $\alpha \to 0, \beta > 0$)**.** *Let $\beta > 0$ be fixed. Assume that $r = r(\alpha)$ is smooth, with $r(0) = r(\alpha = 0)$. Then*

$$\lim_{\alpha \to 0} D_r^{(\alpha,\beta)}[(A+\gamma I), (B+\mu I)] = \frac{1}{\beta^2}\left[\left(\frac{\gamma}{\mu}\right)^{r(0)} - 1\right]\left(1 + r(0)\log\frac{\gamma}{\mu}\right) \quad (59)$$
$$- \frac{1}{\beta^2}\left(\frac{\gamma}{\mu}\right)^{r(0)} \log\det_{2\mathrm{X}}([(B+\mu I)^{-1}(A+\gamma I)]^{r(0)}).$$

Motivated by Theorems 1 and 2, the following is our definition of $D_r^{(\alpha,0)}[(A+\gamma I), (B+\mu I)]$ and $D_r^{(0,\beta)}[(A+\gamma I), (B+\mu I)]$, $\alpha > 0, \beta > 0$.

**Definition 3** (**Limiting cases**)**.** *Let $\alpha, \beta > 0$ be fixed. Let $r \in \mathbb{R}, r \neq 0$ be fixed. For $(A+\gamma I), (B+\mu I) \in \mathscr{PC}_2(\mathcal{H})$, the divergence $D_r^{(\alpha,0)}[(A+\gamma I), (B+\mu I)]$ is defined to be*

$$D_r^{(\alpha,0)}[(A+\gamma I), (B+\mu I)] = \frac{1}{\alpha^2}\left[\left(\frac{\mu}{\gamma}\right)^r - 1\right]\left(1 + r\log\frac{\mu}{\gamma}\right) \quad (60)$$
$$- \frac{1}{\alpha^2}\left(\frac{\mu}{\gamma}\right)^r \log\det_{2\mathrm{X}}([(A+\gamma I)^{-1}(B+\mu I)]^r).$$

*Similarly, the divergence $D_r^{(0,\beta)}[(A+\gamma I), (B+\mu I)]$ is defined to be*

$$D_r^{(0,\beta)}[(A+\gamma I), (B+\mu I)] = \frac{1}{\beta^2}\left[\left(\frac{\gamma}{\mu}\right)^r - 1\right]\left(1 + r\log\frac{\gamma}{\mu}\right) \quad (61)$$
$$- \frac{1}{\beta^2}\left(\frac{\gamma}{\mu}\right)^r \log\det_{2\mathrm{X}}([(B+\mu I)^{-1}(A+\gamma I)]^r).$$

For the case $(A+\gamma I), (B+\mu I) \in \mathscr{PC}_1(\mathcal{H})$, from Definition 3, we recover the formulation stated in Definition 2 in [2], as follows.

**Corollary 2.** *Let $\alpha, \beta > 0$ be fixed. Let $r \in \mathbb{R}, r \neq 0$ be fixed. Assume that $(A+\gamma I), (B+\mu I) \in \mathscr{PC}_1(\mathcal{H})$. Then in Definition 3,*

$$D_r^{(\alpha,0)}[(A+\gamma I), (B+\mu I)] = \frac{r}{\alpha^2}\left[\left(\frac{\mu}{\gamma}\right)^r - 1\right]\log\frac{\mu}{\gamma}$$
$$- \frac{1}{\alpha^2}\left(\frac{\mu}{\gamma}\right)^r \log\det_{\mathrm{X}}[(A+\gamma I)^{-1}(B+\mu I)]^r$$
$$+ \frac{1}{\alpha^2}\mathrm{tr}_X([(A+\gamma I)^{-1}(B+\mu I)]^r - I). \quad (62)$$

$$D_r^{(0,\beta)}[(A+\gamma I),(B+\mu I)] = \frac{r}{\beta^2}\left[\left(\frac{\gamma}{\mu}\right)^r - 1\right]\log\frac{\gamma}{\mu}$$
$$- \frac{1}{\beta^2}\left(\frac{\gamma}{\mu}\right)^r \log\det{}_X[(B+\mu I)^{-1}(A+\gamma I)]^r$$
$$+ \frac{1}{\beta^2}\mathrm{tr}_X([(B+\mu I)^{-1}(A+\gamma I)]^r - I). \quad (63)$$

The following is the generalization of Theorem 9 in [2] to positive definite Hilbert-Schmidt operators.

**Theorem 3** (**Limiting case** $(0,0)$)**.** *Assume that $(A+\gamma I),(B+\mu I) \in \mathscr{PC}_2(\mathcal{H})$. Assume that $r = r(\alpha)$ is smooth, with $r(0) = 0$, $r'(0) \neq 0$, and $r(\alpha) \neq 0$ for $\alpha \neq 0$. Then*

$$\lim_{\alpha \to 0} D_r^{(\alpha,\alpha)}[(A+\gamma I),(B+\mu I)] = \frac{[r'(0)]^2}{8} d_{\mathrm{aiHS}}^2[(A+\gamma I),(B+\mu I)]. \quad (64)$$

*In particular, for $r = 2\alpha$,*

$$\lim_{\alpha \to 0} D_{2\alpha}^{(\alpha,\alpha)}[(A+\gamma I),(B+\mu I)] = \frac{1}{2} d_{\mathrm{aiHS}}^2[(A+\gamma I),(B+\mu I)]. \quad (65)$$

The following is the generalization of Theorem 3 in [2] to positive definite Hilbert-Schmidt operators.

**Theorem 4** (**Symmetric divergences**)**.** *The parametrized family $D_{2\alpha}^{(\alpha,\alpha)}[(A+\gamma I),(B+\mu I)]$, $\alpha \geq 0$, is a family of symmetric divergences on $\mathscr{PC}_2(\mathcal{H})$, with $\alpha = 0$ corresponding to the infinite-dimensional affine-invariant Riemannian distance above and $\alpha = 1/2$ corresponding to the infinite-dimensional symmetric Stein divergence, which is given by $\frac{1}{4}d_{\mathrm{logdet}}^0[(A+\gamma I),(B+\mu I)]$.*

## 6. Properties of the Log-Determinant divergences

The following results establish several important properties of $D_r^{(\alpha,\beta)}$ as defined above, which generalize those from both the finite-dimensional setting [9, 3] and the infinite-dimensional Alpha Log-Det divergences [1] and Alpha-Beta Log-Det divergences [2] for positive definite trace class operators.

In the following theorems, $(A+\gamma I),(B+\mu I) \in \mathscr{PC}_2(\mathcal{H})$.

**Theorem 5** (**Positivity**).
$$D_r^{(\alpha,\beta)}[(A+\gamma I),(B+\mu I)] \geq 0, \tag{66}$$
$$D_r^{(\alpha,\beta)}[(A+\gamma I),(B+\mu I)] = 0 \iff A=B, \gamma=\mu. \tag{67}$$

**Theorem 6** (**Dual symmetry**).
$$D_r^{(\beta,\alpha)}[(B+\mu I),(A+\gamma I)] = D_r^{(\alpha,\beta)}[(A+\gamma I),(B+\mu I)]. \tag{68}$$

*In particular, for $\beta = \alpha$, we have*
$$D_r^{(\alpha,\alpha)}[(B+\mu I),(A+\gamma I)] = D_r^{(\alpha,\alpha)}[(A+\gamma I),(B+\mu I)]. \tag{69}$$

**Theorem 7** (**Dual invariance under inversion**).
$$D_r^{(\alpha,\beta)}[(A+\gamma I)^{-1},(B+\mu I)^{-1}] = D_{-r}^{(\alpha,\beta)}[(A+\gamma I),(B+\mu I)] \tag{70}$$

**Theorem 8** (**Affine invariance**). *For any $(A+\gamma I),(B+\mu I) \in \mathscr{PC}_2(\mathcal{H})$ and any invertible $(C+\nu I) \in \mathrm{HS}_X(\mathcal{H})$, $\nu \neq 0$,*

$$D_r^{(\alpha,\beta)}[(C+\nu I)(A+\gamma I)(C+\nu I)^*, (C+\nu I)(B+\mu I)(C+\nu I)^*]$$
$$= D_r^{(\alpha,\beta)}[(A+\gamma I),(B+\mu I)]. \tag{71}$$

**Theorem 9** (**Invariance under unitary transformations**). *For any $(A+\gamma I),(B+\mu I) \in \mathscr{PC}_2(\mathcal{H})$ and any $C \in \mathcal{L}(\mathcal{H})$, with $CC^* = C^*C = I$,*

$$D_r^{(\alpha,\beta)}[C(A+\gamma I)C^*, C(B+\mu I)C^*] = D_r^{(\alpha,\beta)}[(A+\gamma I),(B+\mu I)]. \tag{72}$$

## 7. Proofs of main results

*7.1. Proofs of the properties of the extended Hilbert-Carleman determinant*

*Proof of Lemma 1*. By definition of the Hilbert-Carleman determinant and the assumption that $I + A > 0$, we have

$$\det_2(I+A) = \det[(I+A)\exp(-A)] = \prod_{k=1}^{\infty}[(1+\lambda_k)\exp(-\lambda_k)] > 0.$$

Thus $\log \det_2(I + A)$ is well-defined and finite, and is given by the series

$$\log \det_2(I + A) = \sum_{k=1}^{\infty} [\log(1 + \lambda_k) - \lambda_k],$$

which necessarily has a finite value.

For the second statement, consider the function $f(x) = \log(1+x) - x$ for $x > -1$. We have $f'(x) = -\frac{x}{1+x}$, with $f'(x) > 0$ for $-1 < x < 0$ and $f'(x) < 0$ for $x > 0$. Thus $f$ has a unique global maximum $f_{\max} = f(0) = 0$. Thus $\forall k \in \mathbb{N}$,

$$\log(1 + \lambda_k) - \lambda_k \leq 0, \quad \text{with equality if and only } \lambda_k = 0.$$

It then follows that $\log(I + A) \leq 0$, with equality if and only $\lambda_k = 0 \; \forall k \in \mathbb{N}$, that is if and only if $A = 0$. $\square$

***Proof of Lemma 2***. We make use of the result that $R_2(A) = (I + A)\exp(-A) - I \in \text{Tr}(\mathcal{H})$ for $A \in \text{HS}(\mathcal{H})$. Thus

$$R_{2,\gamma}(A) = (A + \gamma I)\exp(-A/\gamma) - \gamma I = \gamma\left[\left(\frac{A}{\gamma} + I\right)\exp(-A/\gamma) - I\right] \in \text{Tr}(\mathcal{H}),$$

and hence $R_{2,\gamma}(A) + \gamma I = (A + \gamma I)\exp(-A/\gamma) \in \text{Tr}_X(\mathcal{H})$. Since $R_{2,\gamma}(A) \in \text{Tr}(\mathcal{H})$, the infinite product

$$\prod_{k=1}^{\infty}[(\lambda_k + \gamma)\exp(-\lambda_k/\gamma) - \gamma + 1] = \det[R_{2,\gamma}(A) + I]$$

converges to a finite value. $\square$

***Proof of Lemma 3 (Factorization Rule)***. We have $\left(\frac{A}{\gamma} + I\right)\exp(-A/\gamma) - I \in \text{Tr}(\mathcal{H})$ and thus for the operator

$$(A + \gamma I)\exp(-A/\gamma) = \gamma\left[\left(\frac{A}{\gamma} + I\right)\exp(-A/\gamma) - I\right] + \gamma I \in \text{Tr}_X(\mathcal{H}),$$

its extended Fredholm determinant is given by

$$\det_X[(A + \gamma I)\exp(-A/\gamma)] = \gamma \det\left[\left(\frac{A}{\gamma} + I\right)\exp(-A/\gamma)\right] = \gamma \det_2\left(\frac{A}{\gamma} + I\right).$$

This completes the proof. $\square$

***Proof of Lemma 4 (Similarity Invariant).*** Since $\text{HS}(\mathcal{H})$ is a two-sided ideal in $\mathcal{L}(\mathcal{H})$, we have $CAC^{-1} \in \text{HS}(\mathcal{H})$. Thus

$$C(A+\gamma I)C^{-1} = CAC^{-1} + \gamma I \in \text{HS}_X(\mathcal{H}).$$

By definition of the extended Hilbert-Carleman determinant, we have

$$\det_{2X}[C(A+\gamma I)C^{-1}] = \det_X[C(A+\gamma I)C^{-1}\exp(-CAC^{-1}/\gamma)]$$
$$= \det_X[C(A+\gamma I)C^{-1}(C\exp(-A/\gamma)C^{-1})] = \det_X[C(A+\gamma I)\exp(-A/\gamma)C^{-1}]$$
$$= \det_X[(A+\gamma I)\exp(-A/\gamma)] \quad \text{by Eq. (34)}$$
$$= \det_{2X}(A+\gamma I).$$

This completes the proof. □

***Proof of Lemma 5 (Commutativity).*** Consider the first assumption, that is $(A+\gamma I), (B+\mu I) \in \mathscr{PC}_2(\mathcal{H})$. We write $(A+\gamma I)(B+\mu I)$ and $(B+\mu I)(A+\gamma I)$ as

$$(A+\gamma I)(B+\mu I) = (A+\gamma I)^{1/2}[(A+\gamma I)^{1/2}(B+\mu I)(A+\gamma I)^{1/2}](A+\gamma I)^{-1/2},$$
$$(B+\mu I)(A+\gamma I) = (A+\gamma I)^{-1/2}[(A+\gamma I)^{1/2}(B+\mu I)(A+\gamma I)^{1/2}](A+\gamma I)^{1/2}.$$

By Lemma 4, we then have

$$\det_{2X}[(A+\gamma I)(B+\mu I)] = \det_{2X}[(A+\gamma I)^{1/2}(B+\mu I)(A+\gamma I)^{1/2}]$$
$$= \det_{2X}[(B+\mu I)(A+\gamma I)].$$

The third statement is proved similarly.

Under the second assumption, that is $(A+\gamma I), (B+\mu I) \in \text{Tr}_X(\mathcal{H}), \gamma \neq 0, \mu \neq 0$, we have by definition

$$\det_{2X}[(A+\gamma I)(B+\mu I)] = \det_X[(A+\gamma I)(B+\mu I)]\exp\left(-\frac{\text{tr}[\mu A + \gamma B + AB]}{\gamma\mu}\right)$$
$$= \det_X[(B+\mu I)(A+\gamma I)]\exp\left(-\frac{\text{tr}[\mu A + \gamma B + BA]}{\gamma\mu}\right)$$
$$= \det_{2X}[(B+\mu I)(A+\gamma I)].$$

Here we have made use of the properties $\det_X[(A+\gamma I)(B+\mu I)] = \det_X(A+\gamma I)\det_X(B+\mu I) = \det_X[(B+\mu I)(A+\gamma I)]$ and the commutativity of the trace, namely $\text{tr}(AB) = \text{tr}(BA)$. This completes the proof. □

*Proof of Lemma 6.* We rewrite $(A + \gamma I)(B + \mu I)^{-1}$ as

$$(A + \gamma I)(B + \mu I)^{-1}$$
$$= (B + \mu I)^{1/2}[(B + \mu I)^{-1/2}(A + \gamma I)(B + \mu I)^{-1/2}](B + \mu I)^{-1/2}$$
$$= (B + \mu I)^{1/2}(\Lambda + \nu I)(B + \mu I)^{-1/2}.$$

Similarly,

$$(B + \mu I)^{-1}(A + \gamma I) = (B + \mu I)^{-1/2}(\Lambda + \nu I)(B + \mu I)^{1/2}.$$

By definition of the power function, we then have for any $\alpha \in \mathbb{R}$

$$[(A + \gamma I)(B + \mu I)^{-1}]^\alpha = (B + \mu I)^{1/2}(\Lambda + \nu I)^\alpha (B + \mu I)^{-1/2},$$
$$[(B + \mu I)^{-1}(A + \gamma I)]^\alpha = (B + \mu I)^{-1/2}(\Lambda + \nu I)^\alpha (B + \mu I)^{1/2}.$$

Thus by Lemma 4, we obtain

$$\det_{2\mathrm{X}}([(A + \gamma I)(B + \mu I)^{-1}]^\alpha) = \det_{2\mathrm{X}}[(\Lambda + \nu I)^\alpha]$$
$$= \det_{2\mathrm{X}}([(B + \mu I)^{-1}(A + \gamma I)]^\alpha).$$

This completes the proof. $\square$

**Lemma 9.** *Let $r \neq 0$ be fixed. The function $f(x) = x^r - 1 - r \log(x)$ for $x > 0$ has a unique global minimum $f_{\min} = f(1) = 0$. In other words, $f(x) \geq 0 \; \forall x > 0$, with equality if and only if $x = 1$.*

*Proof of Lemma 9.* We have $f'(x) = \frac{r(x^r - 1)}{x}$. When $r > 0$, we have $x^r < 1$ for $0 < x < 1$ and $x^r > 1$ for $x > 1$. When $r < 0$, we have $x^r > 1$ for $0 < x < 1$ and $x^r < 1$ for $x > 1$. Thus, for all $r \neq 0$, we have $f'(x) < 0$ when $0 < x < 1$ and $f'(x) > 0$ when $x > 1$. Hence $f$ has a unique global minimum $f_{\min} = f(1) = 0$. $\square$

*Proof of Lemma 7.* By Proposition 2 in [15], we have $\log(I + A) \in \mathrm{HS}(\mathcal{H})$ for $(I + A) \in \mathscr{PC}_2(\mathcal{H})$. By definition of the power function, we have

$$(I + A)^\alpha = \exp[\alpha \log(I + A)] = I + \sum_{j=1}^\infty \frac{\alpha^j}{j!}[\log(I + A)]^j.$$

Since HS($\mathcal{H}$) is a Banach algebra under the Hilbert-Schmidt norm, we then have

$$||(I+A)^\alpha - I||_{\text{HS}} = \left\| \sum_{j=1}^\infty \frac{\alpha^j}{j!} [\log(I+A)]^j \right\|_{\text{HS}} \leq \sum_{j=1}^\infty \frac{|\alpha|^j}{j!} ||\log(I+A)||_{\text{HS}}^j$$

$$= \exp(|\alpha| \, ||\log(I+A)||_{\text{HS}}) - 1 < \infty.$$

Thus $(I+A)^\alpha - I \in \text{HS}(\mathcal{H})$. By definition of the Hilbert-Carleman determinant, we then have

$$\det_2[(I+A)^\alpha] = \det[(I+A)^\alpha \exp(-[(I+A)^\alpha - I])]$$

$$= \prod_{k=1}^\infty (1+\lambda_k)^\alpha \exp[1 - (1+\lambda_k)^\alpha] < \infty.$$

Thus the following quantity is well-defined and finite

$$\log \det_2[(I+A)^\alpha] = \sum_{k=1}^\infty [\alpha \log(1+\lambda_k) + 1 - (1+\lambda_k)^\alpha].$$

The statements for the case $I + A \in \mathscr{P}\mathscr{C}_1(\mathcal{H})$ are then obvious from the above series expansions.

By Lemma 9, we have $\forall k \in \mathbb{N}$,

$$\alpha \log(1+\lambda_k) + 1 - (1+\lambda_k)^\alpha \leq 0,$$

with equality if and only if $\lambda_k = 0$. Thus it follows that

$$\log \det_2[(I+A)^\alpha] \leq 0,$$

with equality if and only if $\lambda_k = 0 \; \forall k \in \mathbb{N}$, that is if and only if $A = 0$ (by the assumption that $I + A > 0$). $\square$

*Proof of Lemma 8*. By definition of the power function, we have

$$(A + \gamma I)^\alpha = \exp[\alpha \log(A + \gamma I)] = \exp\left[(\alpha \log \gamma)I + \alpha \log\left(\frac{A}{\gamma} + I\right)\right]$$

$$= \gamma^\alpha \left(\frac{A}{\gamma} + I\right)^\alpha = \gamma^\alpha \left[\left(\frac{A}{\gamma} + I\right)^\alpha - I\right] + \gamma^\alpha I,$$

where $\left[\left(\frac{A}{\gamma} + I\right)^\alpha - I\right] \in \text{Sym}(\mathcal{H}) \cap \text{HS}(\mathcal{H})$ by Lemma 7. Thus it follows that $(A + \gamma I)^\alpha - \gamma^\alpha I \in \text{Sym}(\mathcal{H}) \cap \text{HS}(\mathcal{H})$. Therefore, the extended Hilbert-Carleman determinant of $(A+\gamma I)^\alpha$ is well-defined and finite. By the Factorization Rule (Lemma 3) and Lemma 7, we have

$$\det_{2X}[(A+\gamma I)^\alpha] = \gamma^\alpha \det_2\left[\left(\frac{A}{\gamma}+I\right)^\alpha\right]$$
$$= \gamma^\alpha \prod_{k=1}^\infty \left(1+\frac{\lambda_k}{\gamma}\right)^\alpha \exp\left[1-\left(1+\frac{\lambda_k}{\gamma}\right)^\alpha\right] < \infty.$$

Consequently, the following quantity is also finite

$$\log \det_{2X}[(A+\gamma I)^\alpha] = \alpha \log \gamma + \sum_{k=1}^\infty \left[\alpha \log\left(1+\frac{\lambda_k}{\gamma}\right) + 1 - \left(1+\frac{\lambda_k}{\gamma}\right)^\alpha\right].$$

If $A+\gamma I \in \mathscr{PC}_1(\mathcal{H})$, then $(\frac{A}{\gamma}+I)^\alpha - I \in \text{Tr}(\mathcal{H})$ (see Lemma 6 in [2], or by using a similar argument as in Lemma 7). Thus the following infinite product and series

$$\det\left(\frac{A}{\gamma}+I\right)^\alpha = \prod_{k=1}^\infty \left(\frac{\lambda_k}{\gamma}+1\right)^\alpha,$$
$$\alpha \log \det\left(\frac{A}{\gamma}+I\right) = \alpha \sum_{k=1}^\infty \log\left(\frac{\lambda_k}{\gamma}+1\right),$$
$$\text{tr}\left[\left(\frac{A}{\gamma}+I\right)^\alpha - I\right] = \sum_{k=1}^\infty \left[\left(1+\frac{\lambda_k}{\gamma}\right)^\alpha - 1\right]$$

converge to finite values. These give the last statements of the lemma. $\square$

### 7.2. Proofs for the definition of the Log-Determinant divergences

**Lemma 10.** *Let $\alpha > 0, \beta > 0, \gamma > 0$ be fixed. Let $p, q \in \mathbb{R}$ be such that $p\alpha\gamma^p = q\beta\gamma^{-q}$. Then*

$$\lim_{x \to 0} \frac{1 - \frac{\alpha\gamma^p(1+x)^p + \beta\gamma^{-q}(1+x)^{-q}}{\alpha\gamma^p + \beta\gamma^{-q}}}{x^2} = -\frac{p(p-1)\alpha\gamma^p + q(q+1)\beta\gamma^{-q}}{2(\alpha\gamma^p + \beta\gamma^{-q})}. \quad (73)$$

*Proof of Lemma 10.* Since the limit has the form $\frac{0}{0}$, by L'Hopital's rule, we have

$$\lim_{x \to 0} \frac{1 - \frac{\alpha\gamma^p(1+x)^p + \beta\gamma^{-q}(1+x)^{-q}}{\alpha\gamma^p + \beta\gamma^{-q}}}{x^2}$$
$$= -\frac{1}{2(\alpha\gamma^p + \beta\gamma^{-q})} \lim_{x \to 0} \frac{p\alpha\gamma^p(1+x)^{p-1} - q\beta\gamma^{-q}(1+x)^{-q-1}}{x}.$$

By assumption, we have $p\alpha\gamma^p = q\beta\gamma^{-q}$, so that the previous limit also has the form $\frac{0}{0}$. Applying L'Hopital's rule one more time, we obtain

$$-\frac{1}{2(\alpha\gamma^p + \beta\gamma^{-q})} \lim_{x\to 0}[p(p-1)\alpha\gamma^p(1+x)^{p-2} + q(q+1)\beta\gamma^{-q}(1+x)^{-q-2}]$$
$$= -\frac{p(p-1)\alpha\gamma^p + q(q+1)\beta\gamma^{-q}}{2(\alpha\gamma^p + \beta\gamma^{-q})}.$$

This completes the proof. $\square$

**Corollary 3.** *Let $A \in \text{Sym}(\mathcal{H}) \cap \text{HS}(\mathcal{H})$ be such that $I + A > 0$. Let $\alpha > 0, \beta > 0, \gamma > 0$ be fixed. Let $p, q \in \mathbb{R}$ be such that $p\alpha\gamma^p = q\beta\gamma^{-q}$. Then*

$$I - \frac{\alpha\gamma^p(I+A)^p + \beta\gamma^{-q}(I+A)^{-q}}{\alpha\gamma^p + \beta\gamma^{-q}} \in \text{Sym}(\mathcal{H}) \cap \text{Tr}(\mathcal{H}). \tag{74}$$

*Proof of Corollary 3.* Let $\{\lambda_k\}_{k=1}^\infty$ denote the eigenvalues of $A$ then $\lim_{k\to\infty} \lambda_k = 0$. By Lemma 10, we have

$$\lim_{k\to\infty} \frac{1 - \frac{\alpha\gamma^p(1+\lambda_k)^p + \beta\gamma^{-q}(1+\lambda_k)^{-q}}{\alpha\gamma^p + \beta\gamma^{-q}}}{\lambda_k^2} = -\frac{p(p-1)\alpha\gamma^p + q(q+1)\beta\gamma^{-q}}{2(\alpha\gamma^p + \beta\gamma^{-q})}.$$

This implies that there exists a constant $C > 0$, independent of $k$, and a number $N = N(C) \in \mathbb{N}$, such that

$$\left|1 - \frac{\alpha\gamma^p(1+\lambda_k)^p + \beta\gamma^{-q}(1+\lambda_k)^{-q}}{\alpha\gamma^p + \beta\gamma^{-q}}\right| \leq C\lambda_k^2 \quad \forall k \geq N.$$

Since $\sum_{k=1}^\infty \lambda_k^2 < \infty$ by assumption, it then follows that

$$\sum_{k=1}^\infty \left|1 - \frac{\alpha\gamma^p(1+\lambda_k)^p + \beta\gamma^{-q}(1+\lambda_k)^{-q}}{\alpha\gamma^p + \beta\gamma^{-q}}\right| < \infty,$$

which gives us the desired result. $\square$

*Proof of Proposition 2.* Since $(A + \gamma I), (B + \mu I) \in \mathscr{PC}_2(\mathcal{H})$, we have $\Lambda + \frac{\gamma}{\mu}I = (B + \mu I)^{-1/2}(A + \gamma I)(B + \mu I)^{-1/2} \in \mathscr{PC}_2(\mathcal{H})$, with $\Lambda \in \text{HS}(\mathcal{H})$. Thus it is obvious that $\frac{\alpha(\Lambda + \frac{\gamma}{\mu}I)^p + \beta(\Lambda + \frac{\gamma}{\mu}I)^{-q}}{\alpha + \beta}$ is also positive definite. Let us show that it is an extended trace class operator. Consider the expansion

$$\frac{\alpha(\Lambda + \frac{\gamma}{\mu}I)^p + \beta(\Lambda + \frac{\gamma}{\mu}I)^{-q}}{\alpha + \beta}$$

$$= \frac{\alpha(\frac{\gamma}{\mu})^p + \beta(\frac{\gamma}{\mu})^{-q}}{\alpha + \beta} \frac{\alpha(\frac{\gamma}{\mu})^p(\frac{\mu}{\gamma}\Lambda + I)^p + \beta(\frac{\gamma}{\mu})^{-q}(\frac{\mu}{\gamma}\Lambda + I)^{-q}}{\alpha(\frac{\gamma}{\mu})^p + \beta(\frac{\gamma}{\mu})^{-q}}$$

$$= \frac{\alpha(\frac{\gamma}{\mu})^p + \beta(\frac{\gamma}{\mu})^{-q}}{\alpha + \beta} \left[I - \left(I - \frac{\alpha(\frac{\gamma}{\mu})^p(\frac{\mu}{\gamma}\Lambda + I)^p + \beta(\frac{\gamma}{\mu})^{-q}(\frac{\mu}{\gamma}\Lambda + I)^{-q}}{\alpha(\frac{\gamma}{\mu})^p + \beta(\frac{\gamma}{\mu})^{-q}}\right)\right].$$

By Corollary 3, we have $\left(I - \frac{\alpha(\frac{\gamma}{\mu})^p(\frac{\mu}{\gamma}\Lambda+I)^p + \beta(\frac{\gamma}{\mu})^{-q}(\frac{\mu}{\gamma}\Lambda+I)^{-q}}{\alpha(\frac{\gamma}{\mu})^p + \beta(\frac{\gamma}{\mu})^{-q}}\right) \in \mathrm{Sym}(\mathcal{H}) \cap \mathrm{Tr}(\mathcal{H})$.

Thus it follows that $\frac{\alpha(\Lambda + \frac{\gamma}{\mu}I)^p + \beta(\Lambda + \frac{\gamma}{\mu}I)^{-q}}{\alpha+\beta} \in \mathrm{Tr}_X(\mathcal{H})$. This completes the proof. $\square$

*Proof of Proposition 3*. By Proposition 2, we have

$$\frac{\alpha(\Lambda + \frac{\gamma}{\mu}I)^p + \beta(\Lambda + \frac{\gamma}{\mu}I)^{-q}}{\alpha + \beta} \in \mathscr{PC}_1(\mathcal{H}).$$

Thus its extended Fredholm determinant $\det_X$ is well-defined and finite.

By Proposition 1, we have for any $p \in \mathbb{R}$,

$$[(A + \gamma I)(B + \mu I)^{-1}]^p = (B + \mu I)^{1/2}(\Lambda + \frac{\gamma}{\mu}I)^p(B + \mu I)^{-1/2} \in \mathrm{HS}_X(\mathcal{H}).$$

Thus it follows that

$$\frac{\alpha[(A + \gamma I)(B + \mu I)^{-1}]^p + \beta[(A + \gamma I)(B + \mu I)^{-1}]^{-q}}{\alpha + \beta}$$
$$= (B + \mu I)^{1/2} \left[\frac{\alpha(\Lambda + \frac{\gamma}{\mu}I)^p + \beta(\Lambda + \frac{\gamma}{\mu}I)^{-q}}{\alpha + \beta}\right](B + \mu I)^{-1/2} \in \mathrm{Tr}_X(\mathcal{H}).$$

Thus by Eq. (34), we obtain

$$\det_X \left[\frac{\alpha[(A + \gamma I)(B + \mu I)^{-1}]^p + \beta[(A + \gamma I)(B + \mu I)^{-1}]^{-q}}{\alpha + \beta}\right]$$
$$= \det_X \left[\frac{\alpha(\Lambda + \frac{\gamma}{\mu}I)^p + \beta(\Lambda + \frac{\gamma}{\mu}I)^{-q}}{\alpha + \beta}\right].$$

This completes the proof. $\square$

*Proof of Theorems 1 and 2*. Let $\{\lambda_j\}_{j=1}^\infty$ be the eigenvalues of $\Lambda$. By Theorem 8 in

[2], we have the following expansion

$$D_r^{(\alpha,\beta)}[(A+\gamma I),(B+\mu I)] = \frac{r(\delta - \frac{\alpha}{\alpha+\beta})}{\alpha\beta}\log\left(\frac{\gamma}{\mu}\right) + \frac{1}{\alpha\beta}\log\left(\frac{\alpha(\frac{\gamma}{\mu})^p + \beta(\frac{\gamma}{\mu})^{-q}}{\alpha+\beta}\right)$$

$$+ \frac{1}{\alpha\beta}\log\det\left(\frac{\alpha(\Lambda+\frac{\gamma}{\mu}I)^p + \beta(\Lambda+\frac{\gamma}{\mu}I)^{-q}}{\alpha(\frac{\gamma}{\mu})^p + \beta(\frac{\gamma}{\mu})^{-q}}\right)$$

$$= \frac{r(\delta - \frac{\alpha}{\alpha+\beta})}{\alpha\beta}\log\left(\frac{\gamma}{\mu}\right) + \frac{1}{\alpha\beta}\log\left(\frac{\alpha(\frac{\gamma}{\mu})^p + \beta(\frac{\gamma}{\mu})^{-q}}{\alpha+\beta}\right)$$

$$+ \frac{1}{\alpha\beta}\sum_{j=1}^{\infty}\log\left(\frac{\alpha(\lambda_j+\frac{\gamma}{\mu})^p + \beta(\lambda_j+\frac{\gamma}{\mu})^{-q}}{\alpha(\frac{\gamma}{\mu})^p + \beta(\frac{\gamma}{\mu})^{-q}}\right),$$

where $p = p(\beta) = r(1-\delta) = \frac{r\beta}{\alpha(\frac{\gamma}{\mu})^r + \beta}$, $q = q(\beta) = r\delta = \frac{r\alpha(\frac{\gamma}{\mu})^r}{\alpha(\frac{\gamma}{\mu})^r + \beta}$.

Let $\nu = \frac{\gamma}{\mu}$. By the same argument as in the proof of Theorem 11 in [2], we have

$$\lim_{\beta\to 0} D_r^{(\alpha,\beta)}[(A+\gamma I),(B+\mu I)] = \frac{1}{\alpha^2}[\nu^{-r(0)} + r(0)\log(\nu) - 1]$$

$$+ \frac{1}{\alpha^2}\sum_{j=1}^{\infty}\left[\frac{r(0)}{\nu^{r(0)}}\log\left(\frac{\lambda_j}{\nu}+1\right) + \frac{1}{(\lambda_j+\nu)^{r(0)}} - \frac{1}{\nu^{r(0)}}\right]. \tag{75}$$

By Lemma 8, we have

$$\sum_{j=1}^{\infty}\left[\frac{r(0)}{\nu^{r(0)}}\log\left(\frac{\lambda_j}{\nu}+1\right) + \frac{1}{(\lambda_j+\nu)^{r(0)}} - \frac{1}{\nu^{r(0)}}\right]$$

$$= -\nu^{-r(0)}\sum_{j=1}^{\infty}\left[-r(0)\log\left(\frac{\lambda_j}{\nu}+1\right) - \left(\frac{\lambda_j}{\nu}+1\right)^{-r(0)} + 1\right]$$

$$= -\nu^{-r(0)}(\log\det_{2X}[(\Lambda+\nu I)^{-r(0)}] + r(0)\log\nu).$$

Combining this with Eq. (75), we obtain

$$\lim_{\beta\to 0} D_r^{(\alpha,\beta)}[(A+\gamma I),(B+\mu I)]$$

$$= \frac{1}{\alpha^2}\left[(\nu^{-r(0)}-1)(1 - r(0)\log\nu) - \nu^{-r(0)}\log\det_{2X}[(\Lambda+\nu I)^{-r(0)}]\right].$$

By Lemma 6, we have

$$\det_{2X}[(\Lambda+\nu I)^{-r(0)}] = \det_{2X}([(B+\mu I)^{-1}(A+\gamma I)]^{-r(0)}$$

$$= \det_{2X}([(A+\gamma I)^{-1}(B+\mu I)]^{r(0)}).$$

Substituting this into the previous expression and $\nu = \frac{\gamma}{\mu}$, we obtain the final result.

By dual symmetry, we then obtain $\lim_{\alpha \to 0} D_r^{(\alpha,\beta)}[(A+\gamma I),(B+\mu I)]$ via

$$\lim_{\alpha \to 0} D_r^{(\alpha,\beta)}[(A+\gamma I),(B+\mu I)] = \lim_{\alpha \to 0} D^{(\beta,\alpha)}[(B+\mu I),(A+\gamma I)].$$

This completes the proof. $\square$

*Proof of Corollary 2.* Let us prove the first statement, since the second one is entirely similar. It suffices to prove for $\alpha = 1$. For $(A+\gamma I),(B+\mu I) \in \mathscr{PC}_1(\mathcal{H})$, we have $(\Lambda + \nu I) = (B+\mu I)^{-1/2}(A+\gamma I)(B+\mu I)^{-1/2} \in \mathscr{PC}_1(\mathcal{H})$. By Definition 3,

$$D_r^{(1,0)}[(A+\gamma I),(B+\mu I)] = (\nu^{-r} - 1)(1 - r\log\nu) - \nu^{-r}\log\det{}_{2\mathrm{X}}[(\Lambda + \nu I)^{-r}].$$

By Lemma 8, we have

$$\nu^{-r}\log\det{}_{2\mathrm{X}}[(\Lambda + \nu I)^{-r}] = -r\nu^{-r}\log\det{}_{\mathrm{X}}(\Lambda + \nu I) - \nu^{-r}\mathrm{tr}\left[\left(\frac{\Lambda}{\nu} + I\right)^{-r} - I\right]$$

$$= -r\nu^{-r}\log\det{}_{\mathrm{X}}(\Lambda + \nu I) - \mathrm{tr}[(\Lambda + \nu I)^{-r} - \nu^{-r}I].$$

It then follows that

$$(\nu^{-r} - 1)(1 - r\log\nu) - \nu^{-r}\log\det{}_{2\mathrm{X}}[(\Lambda + \nu I)^{-r}]$$
$$= (\nu^{-r} - 1)(1 - r\log\nu) + r\nu^{-r}\log\det{}_{\mathrm{X}}(\Lambda + \nu I) + \mathrm{tr}[(\Lambda + \nu I)^{-r} - \nu^{-r}I]$$
$$= -r(\nu^{-r} - 1)\log\nu + \nu^{-r}\log\det{}_{\mathrm{X}}(\Lambda + \nu I)^r + (\nu^{-r} - 1 + \mathrm{tr}[(\Lambda + \nu I)^{-r} - \nu^{-r}I])$$
$$= -r(\nu^{-r} - 1)\log\nu - \nu^{-r}\log\det{}_{\mathrm{X}}(\Lambda + \nu I)^{-r} + \mathrm{tr}_X[(\Lambda + \nu I)^{-r} - I].$$

By Lemma 8 in [2], which states that for any $\alpha \in \mathbb{R}$,

$$\det{}_{\mathrm{X}}[(A+\gamma I)(B+\mu I)^{-1}]^\alpha = \det{}_{\mathrm{X}}[(\Lambda + \nu I)^\alpha] = \det{}_{\mathrm{X}}[(B+\mu I)^{-1}(A+\gamma I)]^\alpha,$$
$$\mathrm{tr}_X[(A+\gamma I)(B+\mu I)^{-1}]^\alpha = \mathrm{tr}_X[(\Lambda + \nu I)^\alpha] = \mathrm{tr}_X[(B+\mu I)^{-1}(A+\gamma I)]^\alpha,$$

we have

$$\det{}_{\mathrm{X}}(\Lambda + \nu I)^{-r} = \det{}_{\mathrm{X}}[(B+\mu I)^{-1}(A+\gamma I)]^{-r} = \det{}_{\mathrm{X}}[(A+\gamma I)^{-1}(B+\mu I)]^r,$$
$$\mathrm{tr}_X[(\Lambda + \nu I)^{-r}] = \mathrm{tr}_X[(B+\mu I)^{-1}(A+\gamma I)]^{-r} = \mathrm{tr}_X[(A+\gamma I)^{-1}(B+\mu I)]^r.$$

Combining these with the previous expression, replacing $\nu = \frac{\gamma}{\mu}$, we obtain

$$D_r^{(1,0)}[(A+\gamma I),(B+\mu I)] = \left((\frac{\mu}{\gamma})^r - 1\right)\log\frac{\mu}{\gamma}$$
$$-(\frac{\mu}{\gamma})^r \log\det{}_X[(A+\gamma I)^{-1}(B+\mu I)]^r$$
$$+ \operatorname{tr}_X([(A+\gamma I)^{-1}(B+\mu I)]^r - I).$$

This completes the proof. $\square$

***Proof of Theorem 3.*** The proof is identical to the proof in the setting $(A+\gamma I), (B+\mu I) \in \mathscr{PC}_1(\mathcal{H})$ (Theorem 9 in [2]). $\square$

***Proof of Theorem 4.*** This follows from the dual symmetry in Theorem 6 and the limiting behavior in Theorem 3. $\square$

*7.3. Proofs of the properties of the Log-Determinant divergences*

For the proof on Theorem 5 on positivity, we first need the following technical results.

**Lemma 11.** *Assume that $\gamma > 0, \alpha > 0, \beta > 0$ are fixed. Let $r \in \mathbb{R}, r \neq 0$ be fixed. Then for $\delta = \frac{\alpha\gamma^r}{\alpha\gamma^r+\beta}$, $p = r(1-\delta)$, $q = r\delta$, we have*

$$\frac{r(\delta - \frac{\alpha}{\alpha+\beta})}{\alpha\beta}\log\gamma + \frac{1}{\alpha\beta}\log\left(\frac{\alpha\gamma^p + \beta\gamma^{-q}}{\alpha+\beta}\right) \geq 0. \tag{76}$$

*Equality happens if and only if $\gamma = 1$.*

***Proof of Lemma 11.*** By the strict concavity of the log function, we have

$$\log\left(\frac{\alpha\gamma^p + \beta\gamma^{-q}}{\alpha+\beta}\right) \geq \frac{(p\alpha - q\beta)\log\gamma}{\alpha+\beta},$$

with equality if and only if $\gamma^p = \gamma^{-q} \iff \gamma^{p+q} = \gamma^r = 1$. Since $\gamma > 0$ and $r \neq 0$, this happens if and only if $\gamma = 1$. Thus we have

$$\frac{r(\delta - \frac{\alpha}{\alpha+\beta})}{\alpha\beta}\log\gamma + \frac{1}{\alpha\beta}\log\left(\frac{\alpha\gamma^p + \beta\gamma^{-q}}{\alpha+\beta}\right)$$
$$\geq \frac{1}{\alpha\beta}\left[r(\delta - \frac{\alpha}{\alpha+\beta}) + \frac{p\alpha - q\beta}{\alpha+\beta}\right]\log\gamma$$
$$= \frac{1}{\alpha\beta}\left[q - \frac{(p+q)\alpha}{\alpha+\beta} + \frac{p\alpha - q\beta}{\alpha+\beta}\right]\log\gamma = 0, \text{ since } r = p+q.$$

Equality happens if and only if $\gamma = 1$. $\square$

**Lemma 12.** *Assume that $\gamma > 0, \alpha > 0, \beta > 0$ are fixed. Let $r \in \mathbb{R}, r \neq 0$ be fixed. Assume that $\lambda \in \mathbb{R}$ is also fixed, such that $\lambda + \gamma > 0$. Then for $\delta = \frac{\alpha \gamma^r}{\alpha \gamma^r + \beta}$, $p = r(1 - \delta)$, $q = r\delta$,*

$$\log\left(\frac{\alpha(\lambda + \gamma)^p + \beta(\lambda + \gamma)^{-q}}{\alpha \gamma^p + \beta \gamma^{-q}}\right) \geq 0. \tag{77}$$

*Equality happens if and only if $\lambda = 0$.*

*Proof of Lemma 12.* By the strict concavity of the log function, we have

$$\log\left(\frac{\alpha(\lambda + \gamma)^p + \beta(\lambda + \gamma)^{-q}}{\alpha \gamma^p + \beta \gamma^{-q}}\right)$$

$$= \log\left[\frac{\alpha \gamma^p}{\alpha \gamma^p + \beta \gamma^{-q}}\left(\frac{\lambda}{\gamma} + 1\right)^p + \frac{\beta \gamma^{-q}}{\alpha \gamma^p + \beta \gamma^{-q}}\left(\frac{\lambda}{\gamma} + 1\right)^{-q}\right]$$

$$\geq \frac{p\alpha\gamma^p - q\beta\gamma^{-q}}{\alpha\gamma^p + \beta\gamma^{-q}} \log\left(\frac{\lambda}{\gamma} + 1\right) = 0,$$

since $p\alpha\gamma^p - q\beta\gamma^{-q} = 0$, as can be verified directly using the given hypothesis. Equality happens if and only if $(\frac{\lambda}{\gamma} + 1)^p = (\frac{\lambda}{\gamma} + 1)^{-q} \iff (\frac{\lambda}{\gamma} + 1)^{p+q} = (\frac{\lambda}{\gamma} + 1)^r = 1$. Since $\lambda + \gamma > 0$, $\gamma > 0$, and $r \neq 0$, this happens if and only if $\lambda = 0$. $\square$

*Proof of Theorem 5 (Positivity).* (a) The case $\alpha > 0, \beta > 0$.

Let $\Lambda + \frac{\gamma}{\mu}I = (B + \mu I)^{-1/2}(A + \gamma I)(B + \mu I)^{-1/2}$. Let $\{\lambda_j\}_{j=1}^{\infty}$ be the eigenvalues of $\Lambda$. By Theorem 8 in [2], we have the expansion

$$D_r^{(\alpha,\beta)}[(A + \gamma I), (B + \mu I)] = \frac{r(\delta - \frac{\alpha}{\alpha+\beta})}{\alpha\beta}\log\left(\frac{\gamma}{\mu}\right) + \frac{1}{\alpha\beta}\log\left(\frac{\alpha(\frac{\gamma}{\mu})^p + \beta(\frac{\gamma}{\mu})^{-q}}{\alpha + \beta}\right)$$

$$+ \frac{1}{\alpha\beta}\log\det\left(\frac{\alpha(\Lambda + \frac{\gamma}{\mu}I)^p + \beta(\Lambda + \frac{\gamma}{\mu}I)^{-q}}{\alpha(\frac{\gamma}{\mu})^p + \beta(\frac{\gamma}{\mu})^{-q}}\right)$$

$$= \frac{r(\delta - \frac{\alpha}{\alpha+\beta})}{\alpha\beta}\log\left(\frac{\gamma}{\mu}\right) + \frac{1}{\alpha\beta}\log\left(\frac{\alpha(\frac{\gamma}{\mu})^p + \beta(\frac{\gamma}{\mu})^{-q}}{\alpha + \beta}\right)$$

$$+ \frac{1}{\alpha\beta}\sum_{j=1}^{\infty}\log\left(\frac{\alpha(\lambda_j + \frac{\gamma}{\mu})^p + \beta(\lambda_j + \frac{\gamma}{\mu})^{-q}}{\alpha(\frac{\gamma}{\mu})^p + \beta(\frac{\gamma}{\mu})^{-q}}\right),$$

where $p = p(\beta) = r(1 - \delta) = \frac{r\beta}{\alpha(\frac{\gamma}{\mu})^r + \beta}$, $q = q(\beta) = r\delta = \frac{r\alpha(\frac{\gamma}{\mu})^r}{\alpha(\frac{\gamma}{\mu})^r + \beta}$.

By Lemma 11, we have

$$\frac{r(\delta - \frac{\alpha}{\alpha+\beta})}{\alpha\beta}\log\left(\frac{\gamma}{\mu}\right) + \frac{1}{\alpha\beta}\log\left(\frac{\alpha(\frac{\gamma}{\mu})^p + \beta(\frac{\gamma}{\mu})^{-q}}{\alpha + \beta}\right) \geq 0,$$

with equality if and only if $\frac{\gamma}{\mu} = 1 \iff \gamma = \mu$.

By Lemma 12, we have $\forall j \in \mathbb{N}$,

$$\log\left(\frac{\alpha(\lambda_j + \frac{\gamma}{\mu})^p + \beta(\lambda_j + \frac{\gamma}{\mu})^{-q}}{\alpha(\frac{\gamma}{\mu})^p + \beta(\frac{\gamma}{\mu})^{-q}}\right) \geq 0,$$

with equality if and only $\lambda_j = 0$.

Combining these two results with the previous expression for $D_r^{(\alpha,\beta)}[(A+\gamma I), (B+\mu I)]$, we obtain

$$D_r^{(\alpha,\beta)}[(A + \gamma I), (B + \mu I)] \geq 0,$$

with equality if and only if $\gamma = \mu$ and $\lambda_j = 0 \ \forall j \in \mathbb{N}$, that is $\Lambda = 0$. This is equivalent to $(B + \mu I)^{-1/2}(A + \gamma I)(B + \mu I)^{-1/2} = I \iff (A + \gamma I) = (B + \mu I) \iff A = B, \gamma = \mu$, since $A, B \in \mathrm{HS}(\mathcal{H})$ by assumption.

(b) The case $\alpha = 0, \beta > 0$.

Since the factor $\beta^2$ can be ignored, it suffices to consider the case $\beta = 1$. We have

$$D_r^{(0,1)}[(A + \gamma I), (B + \mu I)] = \left[\left(\frac{\gamma}{\mu}\right)^r - 1\right]\left(1 + r\log\frac{\gamma}{\mu}\right)$$
$$- \left(\frac{\gamma}{\mu}\right)^r \log\det_{2\mathrm{X}}([(B + \mu I)^{-1}(A + \gamma I)]^r).$$

By Lemma 6, we have for any $r \in \mathbb{R}$,

$$\det_{2\mathrm{X}}([(B + \mu I)^{-1}(A + \gamma I)]^r = \det_{2\mathrm{X}}\left[\left(\Lambda + \frac{\gamma}{\mu}I\right)^r\right] = \det_{2\mathrm{X}}\left[\left(\frac{\gamma}{\mu}\right)^r\left(\frac{\mu}{\gamma}\Lambda + I\right)^r\right].$$

By the Factorization Rule in Lemma 3, we then have

$$\det_{2\mathrm{X}}([(B + \mu I)^{-1}(A + \gamma I)]^r = \left(\frac{\gamma}{\mu}\right)^r \det_2\left[\left(\frac{\mu}{\gamma}\Lambda + I\right)^r\right].$$

Combining this with the first expression for $D_r^{(0,1)}[(A + \gamma I), (B + \mu I)]$, we obtain

$$D_r^{(0,1)}[(A + \gamma I), (B + \mu I)] = \left(\frac{\gamma}{\mu}\right)^r - 1 - r\log\frac{\gamma}{\mu}$$
$$- \left(\frac{\gamma}{\mu}\right)^r \log\det_2\left[\left(\frac{\mu}{\gamma}\Lambda + I\right)^r\right].$$

By Lemma 7, we have

$$\log\det_2\left[\left(\frac{\mu}{\gamma}\Lambda + I\right)^r\right] \leq 0, \quad \text{with equality if and only if } \Lambda = 0.$$

By Lemma 9, we have

$$\left(\frac{\gamma}{\mu}\right)^r - 1 - r\log\frac{\gamma}{\mu} \geq 0, \quad \text{with equality if and only if } \frac{\gamma}{\mu} = 1.$$

Together with the previous expression for $D_r^{(0,1)}[(A+\gamma I),(B+\mu I)]$, these imply

$$\det{}_{2\mathrm{X}}([(B+\mu I)^{-1}(A+\gamma I)]^r) \geq 0,$$

with equality if and only if $\Lambda + \frac{\gamma}{\mu}I = (B+\mu I)^{-1/2}(A+\gamma I)(B+\mu I)^{-1/2} = I \iff A+\gamma I = B+\mu I \iff A=B, \gamma=\mu$.

(c) The case $\alpha > 0, \beta = 0$ follows from the previous case by dual symmetry. This completes the proof. □

*Proof of Theorem 6 (Dual symmetry).* For the case $\alpha > 0, \beta > 0$, the proof is identical to that for the setting $(A+\gamma I), (B+\mu I) \in \mathscr{PC}_1(\mathcal{H})$ (Theorem 13 in [2]). The cases $\alpha = 0, \beta > 0$ and $\alpha > 0, \beta = 0$ are obvious from Eqs. (60) and (61). □

*Proof of Theorem 7 (Dual invariance under inversion).* For the case $\alpha > 0, \beta > 0$, the proof is identical to that for the setting $(A+\gamma I), (B+\mu I) \in \mathscr{PC}_1(\mathcal{H})$ (Theorem 14 in [2]).

Consider the case $\alpha = 0, \beta > 0$ (the case $\alpha > 0, \beta = 0$ follows from dual symmetry). It suffices to consider $\beta = 1$. We have

$$(A+\gamma I)^{-1} = \frac{1}{\gamma}I - \frac{A}{\gamma}(A+\gamma I)^{-1}, \quad (B+\mu I)^{-1} = \frac{1}{\mu}I - \frac{B}{\mu}(B+\mu I)^{-1}.$$

By Eq. (61), we have

$$D_r^{(0,1)}[(A+\gamma I)^{-1},(B+\mu I)^{-1}] = \left[\left(\frac{1/\gamma}{1/\mu}\right)^r - 1\right]\left(1 - r\log\frac{1/\gamma}{1/\mu}\right)$$
$$- \left(\frac{1/\gamma}{1/\mu}\right)^r \log\det{}_{2\mathrm{X}}([(B+\mu I)(A+\gamma I)^{-1}]^r)$$
$$= \left[\left(\frac{\mu}{\gamma}\right)^r - 1\right]\left(1 - r\log\frac{\mu}{\gamma}\right) - \left(\frac{\mu}{\gamma}\right)^r \log\det{}_{2\mathrm{X}}([(A+\gamma I)(B+\mu I)^{-1}]^{-r})$$
$$= \left[\left(\frac{\gamma}{\mu}\right)^{-r} - 1\right]\left(1 + r\log\frac{\gamma}{\mu}\right) - \left(\frac{\gamma}{\mu}\right)^{-r} \log\det{}_{2\mathrm{X}}([(B+\mu I)^{-1}(A+\gamma I)]^{-r})$$
$$= D_{-r}^{(0,1)}[(A+\gamma I),(B+\mu I)],$$

where we have used the property $\det_{2X}([(A+\gamma I)(B+\mu I)^{-1}]^{-r}) = \det_{2X}([(B+\mu I)^{-1}(A+\gamma I)]^{-r})$ by Lemma 6. This completes the proof. □

*Proof of Theorem 8 (Affine invariance).* For any $(A+\gamma I), (B+\mu I) \in \mathscr{PC}_2(\mathcal{H})$, and any $(C+\nu I) \in \mathrm{HS}_X(\mathcal{H})$, we have

$$(C+\nu I)(A+\gamma I)(C+\nu I)^*$$
$$= CAC^* + \nu(CA + AC^*) + \nu^2 A + \gamma CC^* + \gamma\nu(C+C^*) + \gamma\nu^2 I \in \mathscr{PC}_2(\mathcal{H}),$$
$$(C+\nu I)(B+\mu I)(C+\nu I)^*$$
$$= CBC^* + \nu(CB + BC^*) + \nu^2 B + \mu CC^* + \mu\nu(C+C^*) + \mu\nu^2 I \in \mathscr{PC}_2(\mathcal{H}),$$

For two operators $(A+\gamma I), (B+\mu I) \in \mathscr{PC}_2(\mathcal{H})$, we then have

$$[(C+\nu I)(A+\gamma I)(C+\nu I)^*][(C+\nu I)(B+\mu I)(C+\nu I)^*]^{-1}$$
$$= (C+\nu I)[(A+\gamma I)(B+\mu I)^{-1}](C+\nu I)^{-1}.$$

Then for any $p \in \mathbb{R}$, we have by Proposition 1,

$$([(C+\nu I)(A+\gamma I)(C+\nu I)^*][(C+\nu I)(B+\mu I)(C+\nu I)^*]^{-1})^p$$
$$= [(C+\nu I)[(A+\gamma I)(B+\mu I)^{-1}](C+\nu I)^{-1}]^p$$
$$= (C+\nu I)[(A+\gamma I)(B+\mu I)^{-1}]^p(C+\nu I)^{-1} \in \mathrm{HS}_X(\mathcal{H})$$

Thus for the cases $\alpha = 0, \beta > 0$ and $\alpha > 0, \beta = 0$, the affine-invariance follows from the Similarity Invariance of the extended Hilbert-Carleman determinant $\det_{2X}$, stated in Lemma 4, along with the invariance of the ratio $\frac{\gamma\nu^2}{\mu\nu^2} = \frac{\gamma}{\mu}$.

For the case $\alpha > 0, \beta > 0$, let $a = \frac{\alpha}{\alpha+\beta}, b = \frac{\beta}{\alpha+\beta}, p = r(1-\delta), q = r\delta$, we have by Proposition 3

$$a[(A+\gamma I)(B+\mu I)^{-1}]^p + b[(A+\gamma I)(B+\mu I)^{-1}]^{-q} \in \mathrm{Tr}_X(\mathcal{H}).$$

It follows then that

$$a([(C+\nu I)(A+\gamma I)(C+\nu I)^*][(C+\nu I)(B+\mu I)(C+\nu I)^*]^{-1})^p$$
$$+ b([(C+\nu I)(A+\gamma I)(C+\nu I)^*][(C+\nu I)(B+\mu I)(C+\nu I)^*]^{-1})^{-q}$$
$$= (C+\nu I)(a[(A+\gamma I)(B+\mu I)^{-1}]^p + b[(A+\gamma I)(B+\mu I)^{-1}]^{-q})(C+\nu I)^{-1}$$
$$\in \mathrm{Tr}_X(\mathcal{H}).$$

From the Similarity Invariance of both the extended Fredholm determinant $\det_X$, stated in Eq. (34), along with the invariance of the ratio $\frac{\gamma\nu^2}{\mu\nu^2} = \frac{\gamma}{\nu}$, we obtain the affine-invariance for $D_r^{(\alpha,\beta)}[(A+\gamma I),(B+\mu I)]$. □

*Proof of Theorem 9 (Unitary invariance).* The proof for this theorem is similar to that of Theorem 9, by utilizing the fact that $C^* = C^{-1}$ and the Similarity Invariance

$$\det_X[C(A+\gamma I)C^{-1}] = \det_X(A+\gamma I), \quad A+\gamma I \in \text{Tr}_X(\mathcal{H}),$$

for the case $\alpha > 0, \beta > 0$, and

$$\det_{2X}[C(A+\gamma I)C^{-1}] = \det_{2X}(A+\gamma I), \quad A+\gamma I \in \text{HS}_X(\mathcal{H}),$$

for the cases $\alpha > 0, \beta = 0$ and $\alpha = 0, \beta > 0$. □